\documentclass[12pt]{amsart}
\usepackage{latexsym, amsmath, amscd, amssymb, amsthm} 

\newtheorem{te}{Theorem}[section]

 \newtheorem{lm}{Lema}[section] 
 
\begin{document}

\noindent
\title{Bivariate Poincar\'e series for  algebra of covariants   of a binary form }

\author{Leonid  Bedratyuk}

\address{ Khmelnytsky National University, Instituts'ka st. 11, Khmelnytsky, 29016, Ukraine}

\begin{abstract}
A  formula for computation of the bivariate Poincar\'e series $\mathcal{P}_d(z,t)$ for  the algebra of  covariants of binary $d$-form is found. 
\end{abstract}

\maketitle

\section{Introduction}
\label{1}
Let $V_d$ be the complex vector space  of  binary forms of degree $d$ endowed with the natural action of the special linear group  $G=SL(2,\mathbb{C}).$  Consider the corresponding action of the group  $G$ on  the coordinate rings $\mathbb{C}[V_d]$ and $\mathbb{C}[V_d \oplus \mathbb{C}^2 ].$ 
Denote   by ${\mathcal{I}_d=\mathbb{C}[V_d  ]^{\,G}}$ and by ${\mathcal{C}_d=\mathbb{C}[V_d \oplus \mathbb{C}^2 ]^{\,G}}$ the  subalgebras of   $G$-invariant polynomial functions. 
In the language  of classical invariant theory  the algebras  $\mathcal{I}_d$ and   $\mathcal{C}_d$ are called the algebra  of invariants and the algebra of covariants for  the  binary form of  degree $d,$ respectively.
 The algebra  $\mathcal{C}_{d}$ is a finitely  generated bigraded algebra: 
$$
\mathcal{C}_{d}=(\mathcal{C}_{d})_{0,0}+(\mathcal{C}_{d})_{1,0}+\cdots+(\mathcal{C}_{d})_{i,j}+ \cdots,
$$
where  each subspace   $(\mathcal{C}_{d})_{i,j}$ of covariants of degree $i$  and order $j$ is   finite-dimensional. The formal power series 
 $\mathcal{P}_{d}(z,t) \in \mathbb{Z}[[z,t]],$
$$
\mathcal{P}_{d}(z,t)=\sum_{i,j=0}^{\infty }\dim((\mathcal{C}_{d})_{i,j}) z^i t^j,
$$ 
is called the bivariate Poincar\'e series   of the algebra of   covariants  $\mathcal{C}_{d}.$
It is clear that the series  $\mathcal{P}_{d}(z,0)$ is the Poincare series of the algebra $\mathcal{I}_d$  and the series $\mathcal{P}_{d}(z,1)$ is the Poincare series of the algebra $\mathcal{C}_d$   with respect to the usual grading of the algebras by degree.
 The finitely  generation of the algebra of   covariants implies  that its bivariate Poincar\'e series is the  power series expansion  of a  rational function.  We consider here the problem of
computing efficiently this rational function.

 Calculating  the  Poincar\'e series   of the  algebras of invariants and covariants was an  important  object of research  in  invariant theory in  the 19th century.
For  the cases  $d\leq 10,$ $d=12$ the  series $\mathcal{P}_{d}(z,t)$  were calculated by Sylvester, see  in \cite{SF, Sylv-12} the big tables of $\mathcal{P}_{d}(z,t),$ named them as Generating Functions for covariants, reduced form. All those calculations are correct  up to $d=6.$

Relatively recently, Springer \cite{SP} found an  explicit formula for computing the Poincar\'e  series of the algebra of invariants $\mathcal{I}_d.$ 
 In the paper we  have  proved a  Cayley-Sylvester type formula  for calculating of $\dim (\mathcal{C}_{d})_{i,j}$   and a Springer type formula  for calculation of $\mathcal{P}_d(z,t).$ By using the  formula, the series $\mathcal{P}_d(z,t)$  is calculated  for $d\leq 20.$

\section{ Cayley-Sylvester type formula  for  $\dim (\mathcal{C}_{d})_{i,j}$}

 To begin, we give a proof of the Cayley-Sylvester type  formula for the dimension if the  graded component $(\mathcal{C}_{d})_{i,j}.$

Let  $V_d=\langle v_0,v_1,...,v_d \rangle,$ $\dim V_d=d+1$ be standard  irreducible  representation of the Lie algebra  $\mathfrak{sl_{2}}.$
The basis elements    $ \left( \begin{array}{ll}  0\, 1 \\ 0\,0 \end{array} \right),$ $ \left( \begin{array}{ll}  0\, 0 \\ 1\,0 \end{array} \right)$, $ \left( \begin{array}{ll}  1 &  \phantom{-}0 \\  0 &-1 \end{array} \right)$ of the algebra    $\mathfrak{sl_{2}}$ act on    $V_d$  by the derivations  $D_1, D_2, E$ : 
$$
D_1(v_i)=i\, v_{i-1},  D_2(v_i)=(d-i)\,v_{i+1}, E(v_i)=(d-2\,i)\,v_i.
$$
The action of   $\mathfrak{sl_{2}}$  is extended to an action on the symmetrical algebra  $S(V_d)$ in the natural way. 

Let   $\mathfrak{u}_{2}$ be  the maximal unipotent subalgebra of $\mathfrak{sl}_{2}.$ The algebra  $\mathcal{S}_{d},$ defined by 
$$
\mathcal{S}_{d}:= \displaystyle{ S(V_d)^{\mathfrak{u_{2}}}}=\{ v \in S(V_d)|  D_1(v)=0 \},
$$
is called the {\it  algebra of semi-invariants}  of the binary form  of degree $d.$ For any element $v \in \mathcal{S}_{d}$ a natural number $s$ is called {\it the  order} of the element $v$ if the number $s$ is the smallest natural number such that \begin{equation*}D_2^s(v) \ne 0, D_2^{s+1}(v) = 0.\end{equation*}
It is  clear that any semi-invariant  $v \in \mathcal{S}_d$  of  order $i$ is the highest weight vector  for an  irreducible $\mathfrak{sl_{2}}$-module   of the dimension $i+1$ in $S(V_d).$

The classical  theorem \cite{Rob} of Roberts implies  an isomorphism of the algebra of covariants and the algebra of semi-invariants. Thus, it is  enough to compute the Poincar\'e  series of the algebra $\mathcal{S}_d.$ 

The algebra  $S(V_d)$  is $\mathbb{N}$-graded 
$$
S(V_d)=S^0(V_d)+S^1(V_d)+\cdots +S^n(V_d)+\cdots,
$$
and each   $S^n(V_d)$ is a completely  reducible
 representation of the Lie algebra  $\mathfrak{sl_{2}}.$
Thus, the following decomposition  holds
$$
S^n(V_d) \cong \gamma_d(n,0) V_0+\gamma_d(n,1) V_1+ \cdots +\gamma_d(n,d\cdot n) V_{d\cdot n},  \eqno{(*)}
$$
here  $\gamma_d(n,k)$ is  the  multiplicity of the representation  $V_k$  in the decomposition of  $S^n(V_d).$ On the other hand, the multiplicity  $\gamma_d(i,j)$  of the  representation  $V_j$ is  equal to the number of linearly independent homogeneous semi-invariants of  degree $i$   and order $j$  for the  binary $d$-form.
This argument proves 
\begin{lm}
$$
\dim (\mathcal{C}_d)_{i,j}=\gamma_d(i,j).
$$
\end{lm}
The set of weights (eigenvalues of the operator $E$) of a representation  $W$ denote by  $\Lambda_{W},$  in particular, $\Lambda_{V_d}=\{-d, -d+2,  \ldots, d-2, d \}.$ 

A  formal sum 
$$
{\rm Char}(W)=\sum_{k \in \Lambda_{W}} n_W(k) q^k,
$$
is called the character   of a representation  $W,$  
here   $n_W(k)$ denotes the   multiplicity  of the weight $k \in \Lambda_{W}.$
Since, the    multiplicity of any weight of the irreducible representation $V_d$  is  equal to 1, we have  
$$
{\rm Char}(V_d)=q^{-d}+q^{-d+2}+\cdots+q^{d-2}+q^{d}.
$$

The  character   $ {\rm Char}(S^n(V_d))$ of the representation  $S^n(V_d)$  equals   $$H_n(q^{-d},q^{-d+2},\ldots,q^{d}),$$ (see \cite{FH}),  where   $H_n(x_0,x_1,\ldots,x_d)$ is  the complete symmetrical function    $$H_n(x_0,x_1,\ldots,x_d)=\sum_{|\alpha|=n} x_0^{\alpha_0} x_1^{\alpha_1} \ldots x_d^{\alpha_d} , |\alpha|=\sum_i \alpha_i.$$

By replacing    $x_k$ with $q^{d-2\,k},$ $k=0,\ldots, d,$  we  obtain the specialized expression for   ${\rm Char}(S^n(V_d)):$ 
$$
{\rm Char}(S^n(V_d))= \sum_{|\alpha|=n} (q^d)^{\alpha_0} (q^{d-2\cdot 1})^{\alpha_1} \ldots (q^{d-2\,d})^{\alpha_d} =
$$
$$
= \sum_{|\alpha|=n} q^{d\,n-2 (\alpha_1+2\alpha_2+\cdots + d\, \alpha_d)}=\sum_{k=0}^{d\,n} \omega_d(n,k) q^{d\,n-2\,k},
$$
here   $\omega_d(n,k) $  is the number non-negative integer solutions of the equation $$\alpha_1+2\alpha_2+\cdots + d\, \alpha_d=\displaystyle \frac{d\,n-k}{2}$$  on the assumption that  $ \alpha_0+\alpha_1+\cdots +  \alpha_d=n.$ In particular, the coefficient of $q^0$ (the     multiplicity  of zero weight ) is  equal to  $ \displaystyle \omega_d \left(n,\frac{d\,n}{2}\right),$ and the coefficient of  $q^1$ is  equal to   $\displaystyle  \omega_d \left(n,\frac{d\,n-1}{2}\right).$

On the other hand, the decomposition $(*)$ implies the  equality for the characters:  
$$
{\rm Char}(S^n(V_d))=\gamma_d(n,0) {\rm Char}(V_0) +\gamma_d(n,1) {\rm Char}(V_1)+ \cdots +\gamma_d(n, d\,n) {\rm Char}(V_{d\,n}).
$$
We can summarize what we have shown so far in 
\begin{te} $$\dim (\mathcal{C}_d)_{i,j}= \omega_d \left(i,\frac{d\,i-j}{2}\right)- \omega_d \left(i,\frac{d\,i-(j+2)}{2}\right)
.$$
\end{te}
\begin{proof}
The  weight $j$  appears  once  in any  representation $V_k,$  for  $k=j \mod 2,$ $k\geq$  therefore 
$$
 \omega_d \left(i,\frac{d\,i-j}{2}\right)=\gamma_d(i,j)  +\gamma_d(i,j+2) + \cdots +\gamma_d(i, j+4)+...
$$
Similarly, 
$$
 \omega_d \left(i,\frac{d\,i-(j+2)}{2}\right)=\gamma_d(i,j+2)  +\gamma_d(i,j+4) + \cdots +\gamma_d(i, j+6)+...
$$
Thus,
 \begin{gather*}
 \omega_d \left(i,\frac{d\,i-j}{2}\right)- \omega_d \left(i,\frac{d\,i-(j+2)}{2}\right)=\gamma_d(i,j).
\end{gather*}
By using Lemma 1 we obtain 
$$\dim (\mathcal{C}_d)_n= \omega_d \left(n,\frac{d\,n}{2}\right)+ \omega_d \left(n,\frac{d\,n-1}{2}\right)
.$$

\end{proof}

Note,  that the original Cayley-Sylvester formula is 
$$
 \dim (\mathcal{I}_d)_n= \omega_d \left(n,\frac{d\,n}{2}\right)- \omega_d \left(n,\frac{d\,n}{2}-1\right).
$$
Also,  in \cite{BC1}, we  proved that 
$$
 \dim (\mathcal{C}_d)_n= \omega_d \left(n,\frac{d\,n}{2}\right)+ \omega_d \left(n,\frac{d\,n-1}{2}\right).
$$

\section{Calculation  of $\dim (C_{d})_{i,j}$}

 It  is  well known  that  the number
 $\displaystyle \omega_d \left(i,\frac{d\,i-j}{2}\right)$  of   non-negative integer solutions of the following system
$$
\left \{
\begin{array}{c}
\displaystyle \alpha_1+2\alpha_2+\cdots + d\, \alpha_d=\displaystyle \frac{d\,i-j}{2},\\
\\
\displaystyle \alpha_0+\alpha_1+\cdots + \alpha_d=i
\end{array}
\right.
$$
is given by the  coefficient of $\displaystyle z^n t^{\frac{d\,i-j}{2}} $  of the  generating function

$$
f_{d}(z,t)=\frac{1}{(1-z )(1-z\,t)\ldots (1-z\,t^d)}.
$$
We will use the notation $[x^k] F(x)$ to denote the coefficient of $x^k$  in the series expansion of $F(x) \in \mathbb{C}[[x]].$ Thus 
  $$\omega_d \left(i,\frac{d\,i-j}{2}\right)=\left[  z^i t^{\frac{d\,i-j}{2}}\right]f_{d}(z,t).$$
  It is  clear that 
$$\omega_d \left(i,\frac{d\,i-j}{2}\right)=\left[  z^i t^{d\,i-j}\right]f_{d}(z,t^2)=\left[ (z\,t^d)^n\right] t^j f_{d}(z,t^2).$$

Similarly, the number  $\omega_d \left(i,\frac{d\,i-(j+2)}{2}\right)$ of   non-negative integer solutions of the following system
 
$$
\left \{
\begin{array}{c}
\displaystyle \alpha_1+2\alpha_2+\cdots + d\, \alpha_d=\displaystyle \frac{d\,i-(j+2)}{2},\\
\\
\displaystyle \alpha_0+\alpha_1+\cdots + \alpha_d=i
\end{array}
\right.
$$
equals  $$\left[z^i t^{\frac{d\,i-(j+2)}{2}}\right]f_{d}(z,t)=\left[ z^i t^{d\,i-(j+2)}\right] f_{d}(z,t^2)=\left[ (z t^d)^i\right]t^{j+2}\,f_{d}(z,t^2).$$
Therefore,
\begin{gather*}
\omega_d \left(i,\frac{d\,i-j}{2}\right)- \omega_d \left(i,\frac{d\,i-(j+2)}{2}\right)=\\=\left[ (z\,t^d)^n\right] t^j f_{d}(z,t^2)-\left[ (z\,t^d)^i\right] t^{j+2} f_{d}(z,t^2)=\\ \\
=\left[ (z\,t^d)^i\right] (t^j-t^{j+2}) f_{d}(z,t^2)=\left[ (z\,t^d)^i\right] t^j(1-t^2) f_{d}(z,t^2)=\\ \\=\left[ z^i t^{d\,i-j}\right](1-t^2) f_{d}(z,t^2)
\end{gather*}
Thus,  the  following statement holds.
\begin{te}  The number  $\dim (\mathcal{C}_d)_{i,j}$ of   linearly independent  covariants of degree $i$ and order $j$ for the  binary  $d$-form is given  by the formula
\begin{gather*}
\dim (C_d)_{i,j}=\left[ z^i t^{d\,i-j}\right]\left( \frac{1-t^2}{(1-z )(1-z\,t^2)\ldots (1-z\,t^{2\,d})}\right)=\\
=\left[ z^i t^{\frac{d\,i-j}{2}}\right]\left( \frac{1-t}{(1-z )(1-z\,t)\ldots (1-z\,t^{d})}\right).
\end{gather*}
\end{te}
By using the decomposition 
$$
 \frac{1}{(1-z )(1-z\,t)\ldots (1-z\,t^{d})}=\sum_{k=0}^{\infty} \left[ \begin{array}{c} d \\ i  \end{array} \right]_t z^i,
$$
where $\left[ \begin{array}{c} d \\ n  \end{array} \right]_q$ is the  $q$-binomial   coefficient:
$$
\left[ \begin{array}{c} d \\ n  \end{array} \right]_q:=\frac{(1-q^{d+1})(1-q^{d+2})\ldots(1-q^{d+n})}{(1-q)(1-q^{2})\ldots(1-q^{n})},
$$
we obtain  the well-known formula
$$
\dim (C_d)_{i,j}=\left[  t^{\frac{d\,i-j}{2}}\right] (1-t) \left[ \begin{array}{c} d \\ i  \end{array} \right]_t,
$$
for  instance, see \cite{SP}.

\section{ Explicit formula  for $\mathcal{P}_d(z,t)$}

 Let us prove  Springer-type  formula for the Poincar\'e  series  $\mathcal{P}_d(z,t)$ of the  algebra covariants of  the binary $d$-form.  Consider the $\mathbb{C}$-algebra $\mathbb{Z}[[t,z]]$   of  formal   power series.
For an integer    $n \in \mathbb{N}$  define the  $\mathbb{C}$-linear function
$$ \Psi_{n}: \mathbb{Z}[[t,z]] \to \mathbb{Z}[[t,z]],$$  in the  following  way
$$
\Psi_{n}\bigl(z^{i} t^{j}\bigr)=\left \{ \begin{array}{l} z^i t^{n\, i- j}, \text{  if  } \displaystyle  n \,i-j \geq 0,\\ 0, \text{ otherwise. } \end{array} \right. 
$$

The main idea of these calculations is that 
the  Poincar\'e series $ \mathcal{P}_ {d} (z,t) $ can be expressed  in terms of functions $ \Psi.$ The following simple but important statement  holds.
\begin{lm} 
$$
\mathcal{P}_{d}(z,t)=\Psi_{d}\left(\frac{1-t^2}{(1-z\,t)(1-z t^2)\ldots (1-z t^{2d})}\right).
$$
\end{lm}
\begin{proof} Theorem 2   implies  that   $\dim(C_{d})_{i,j}=[z^i t^{di-j}]f_d(z,t^2).$ 
Then
$$
\begin{array}{l}
\displaystyle \mathcal{P}_{d}(z,t) = \sum_{i,j=0}^{\infty}  \dim(C_{d})_{i,j} z^i t^j=\sum_{i,j=0}^{\infty} \bigl([(z^i t^{d\,i-j})]f_d(z,t^2)\bigr)z^i t^j{=}\\ \\
= \Psi_{d}(f_d(z,t^2)).
\end{array}
$$
\end{proof}

Let  $\psi_n:\mathbb{Z}[[t]] \to \mathbb{Z}[[t,z]]$ be  a  $\mathbb{C}$-linear function   defined by  
$$
\psi_{n}\left(t^m\right):=\left \{ \begin{array}{l} z^i t^j, \text{  if  } m=n\,i-j, j<n,  \\ 0, \text{ otherwise }  \end{array} \right. 
$$
Note that $\psi_n(t^0)=1$ and $\psi_1(t^m)=z^m,$ $\psi_0(t^m)=1,$  Also put $\psi_{n}\left(t^m\right)=0$  for $n<0.$

It is clear that for arbitrary series   we  obtain
$$
\psi_{n}\left( \sum_{m=0}^{\infty} a_m t^m \right)=\psi_{n}\left( \sum_{n\,i-j\geq 0}^{\infty} a_{n\,i-j} t^{n\,i-j} \right)=\sum_{i,j=0}^{\infty} a_{n\,i-j}z^i t^j.
$$

In  some  cases   calculating  of the functions $ \Psi$  can  be reduced   to  calculating of  the functions    $\psi$.  The following  statements hold: 
 
\begin{lm}\label{lm3}
\begin{align*}
(i) & \text{ For }  R(t), H(t)  \in \mathbb{C}[[t]]  \text{    holds  } \psi_n( R(t^n) H(t))=R(z) \psi_n(H(t));
\\
(ii) & \text{ for }  R(t) \in \mathbb{C}[[t]]   \text{    and for }  n, k  \in \mathbb{N}     \text{  holds} 
\\
&  \displaystyle  \Psi_{n}\left( \frac{R(t)}{1-z t^k} \right)=\left \{ \begin{array}{l} \displaystyle \frac{\psi_{n-k}(R(t))}{1-zt^{n-k}}, n \geq k, \\   \\ 0,  \text{ if  } n < k. \end{array} \right. \end{align*}
\end{lm}
\begin{proof} 
 \noindent
$(i)$ The  statement follows from the  linearity of the function $\psi_n$ and from the following simple observation 
$$
\psi_n( t^{nk} t^{ni-j})=\psi_n( t^{n(k+i)-j})=z^{k+i} t^j= z^k z^i t^j =z^k \psi_n(t^{ni-j}).
$$

 \noindent
$(ii)$
Let  $ \displaystyle R(z)=\sum_{n\,i-j \geq 0} a_{n\,i-j} t^{n\,i-j}.$  Then for   $k < n$  we  have 
$$
\begin{array}{l}
\displaystyle \Psi_{n}\left( \frac{R(z)}{1-z t^k} \right)=\Psi_{n}\Big( \sum_{n\,i-j \geq 0} \sum_{s=0}^{\infty} a_{n\,i-j}  t^{n\,i-j} (z t^k)^s\Big)=\\
\\
\displaystyle =\Psi_{n}\Big(\sum_{n\,i-j \geq 0} \sum_{s=0}^{\infty} a_{n\,i-j} z^s t^{n\,i-j+ks} \Big){=}\sum_{n\,i-j \geq 0} \sum_{s=0}^{\infty} a_{n\,i-j} z^i \Psi_n \left(z^s t^{ks-j})\right)=\\ \\
=\displaystyle \sum_{n\,i-j \geq 0} \sum_{s=0}^{\infty} a_{n\,i-j} z^{i+s} t^{(n-k)s+j}=\sum_{n\,i-j \geq 0} \sum_{s=0}^{\infty} a_{n\,i-j}  z^i t^j (z t^{(n-k)})^s =\\ \\
=\displaystyle  \sum_{n\,i-j \geq 0}a_{n\,i-j}z^i t^j \frac{1}{1-z t^{n-k}}=\frac{\psi_{n-k}(R(t))}{1-zt^{n-k}}.
\end{array}
$$
 \end{proof}

Now we  can present    Springer type  formula  for calculating of the bivariate Poincar\'e  series $\mathcal{P}_{d}(z,t)$
\begin{te}
\begin{align*}
\displaystyle \mathcal{P}_{d}(z,t)=\sum_{0\leq k <d/2} \psi_{d-2\,k} \left( \frac{(-1)^k t^{k(k+1)} (1-t^2)}{(t^2,t^2)_k\,(t^2,t^2)_{d-k}} \right) \frac{1}{1-z t^{d-2\,k}},
\end{align*}
here $ (a,q)_n=(1-a) (1-a\,q)\cdots (1-a\,q^{n-1}) $  is $q$-shifted factorial.
\end{te}
\begin{proof}
Consider the partial fraction decomposition of  the rational function $f_d(z,t^2):$
$$
f_d(z,t^2)=\sum_{k=0}^{d} \frac{R_k(z)}{1-t z^{2\,k}}.
$$
It  is easy  to  see,  that 
$$
\begin{array}{l}
\displaystyle  R_k(t)=\lim_{z \to t^{-2k}}\left( f_d(z,t^2)(1-z t^k) \right)=\lim_{z \to t^{-2\,k}}\left( \frac{(1-t^2)}{(z,t)_{d+1}}(1-z t^{2\,k}) \right)=\\
\\
\displaystyle =\frac{1-t^2}{(1-t^{-2k})(1-t^{2-2k})\cdots (1-t^{2(k-1)-2k})(1-t^{2(k+1)-2k}) \cdots (1-t^{2d-2k})}=\\
\\
\displaystyle =\frac{t^{2k+(2k-2)+\ldots+2}(1-t^2)}{(t^{2k}-1)(t^{2k-2}-1)\cdots (t^2-1)(1-t^2) \cdots (1-t^{2d-2k}) }=\\
\\
\displaystyle=\frac{(-1)^k t^{k(k+1)}(1-t^2)}{(t^2,t^2)_{k} (t^2,t^2)_{d-k}}.
\end{array}
$$
Using  the above Lemmas we obtain
$$
\begin{array}{l}
\displaystyle \mathcal{P}_{d}(z,t) {=}\Psi_{d}\bigl(f_d(z,t^2)\bigr){=}\Psi_{d}\left( \sum_{k=0}^{n} \frac{R_k(t^2)}{1-z t^{2\,k}} \right){=} \\
\\
\displaystyle =\sum_{0\leq k <d/2} \varphi_{d-2\,k} \left( \frac{(-1)^k t^{k(k+1)}(1-t^2)}{(t^2,t^2)_k (t^2,t^2)_{d-k}} \right) \frac{1}{1-z t^{d-2k}}.
\end{array}
$$
\end{proof}

  For direct computations we use the following technical Lemma

\begin{lm} For    $R(t) \in \mathbb{C}[[t]]$  we  have 
$$
\psi_n\left(\frac{R(t)}{(1-t^{k_1})(1-t^{k_2})\cdots(1-t^{k_m})} \right)= 
\frac{\psi_n\bigr(R(z)Q_n(t^{k_1})Q_n(t^{k_2})Q_n(t^{k_m})\bigr)}{(1-z^{k_1})(1-z^{k_2})\cdots(1-z^{k_m})},
$$
here  $Q_n(t)=1+t+t^2+\ldots+t^{n-1},$  and  $k_i$ are natural numbers.
\end{lm}
\begin{proof}
Taking into account the Lemma \ref{lm3}  we get 
$$
\psi_n\left(\frac{g(t)}{1-t^{m}} \right)=\psi_n\left(\frac{g(t)}{1-t^{n\,m}} \frac{1-t^{nm}}{1-t^m} \right)=\frac{1}{1-t^m}\psi_n \left( g(t) \frac{1-t^{nm}}{1-t^m} \right)=
$$
\begin{multline*}
=\frac{1}{1-t^m}\psi_n \left( g(t)(1+t^m+(t^m)^2+\cdots+(t^m)^{n-1}) \right)=
\frac{1}{1-t^m}\psi_n \left( g(t)Q_n(t^m) \right).
\end{multline*}
In a similar fashion  we prove the general case.
\end{proof}

By using Lemma  $4$ the bivariate Poincar\'e series  $\mathcal{P}_d(z,t)$   for  $d \leq 20$ are   found. All these  results agree  with Sylvester's calculations  up to $d=6,$  see \cite{SF}, \cite{Sylv-12}.

 Below  list of several series:

$$
\begin{array}{l}
\mathcal{P}_1(z,t) = {\displaystyle \frac {1}{1-z t}},
\mathcal{P}_2(z)= \displaystyle \frac {1}{(1-z t^2) \, (1-z^2)},\\
\displaystyle \mathcal{P}_{3}(z,t) :={\frac {{z}^{2}{t}^{2}-zt+1}{ \left( 1-zt \right)  \left( 1-z{t}^{3}
 \right)  \left( 1-{z}^{4} \right) }},\\
\displaystyle \mathcal{P}_{4}(z,t)  ={\frac {{z}^{2}{t}^{4}-z{t}^{2}+1}{ \left( 1-z{t}^{2} \right) 
 \left( 1-{t}^{4}z \right)  \left( 1-{z}^{2} \right)  \left( 1-{z}^
{3} \right) }},\\
\displaystyle \mathcal{P}_{5}(z,t)={\frac {p_5(z,t)}{ \left( 1-zt \right)  \left( 1-{t}^{3}z
 \right)  \left( 1-z{t}^{5} \right)  \left( 1-{z}^{8} \right) 
 \left( 1-{z}^{6} \right)  \left( 1-{z}^{4} \right) }}
\end{array}
$$
\begin{gather*}
p_5(z,t)=1+{z}^{7}{t}^{3}-{z}^{6}{t}^{4}+{z}^{2}{t}^{2}+2\,{z}^{7}t-{z}
^{5}{t}^{5}-{z}^{8}{t}^{2}-2\,{z}^{8}{t}^{6}-{z}^{8}{t}^{4}+{z}^{5}{t}
^{3}+\\+{z}^{5}t+{z}^{9}{t}^{7}-{z}^{10}{t}^{6}+{z}^{10}{t}^{2}-{z}^{10}{
t}^{4}-{z}^{11}{t}^{3}+{z}^{9}{t}^{3}-{t}^{3}z-{z}^{6}+{z}^{4}{t}^{4}-
zt+\\+{z}^{2}{t}^{6}+{z}^{2}{t}^{4}+{z}^{12}+{z}^{14}{t}^{6}-{z}^{13}t-{z
}^{13}{t}^{5}-{z}^{13}{t}^{3}-{z}^{15}{t}^{7}+{z}^{14}{t}^{4}-{z}^{3}{
t}^{7}+{z}^{7}{t}^{5}.
\end{gather*}

\end{document}